\documentclass{article}

\usepackage{arxiv}

\usepackage[utf8]{inputenc} 
\usepackage[T1]{fontenc}    
\usepackage{hyperref}       
\usepackage{url}            
\usepackage{booktabs}       
\usepackage{amsfonts}       
\usepackage{nicefrac}       
\usepackage{microtype}      
\usepackage{lipsum}
\usepackage{graphicx}
\usepackage{amsmath}
\graphicspath{ {./images/} }

\title{Two-Dimensional Quaternion  Linear Canonical Transform: A Novel Framework for Probability Modeling}

\author{
 Muhammad Adnan Samad \\
  School of Automation\\
  Beijing Institute of Technology\\
  Beijing, China \\
  \texttt{adnan.samad@bit.edu.cn} \\
   \And
Yuanqing Xia \\
  School of Automation\\
  Beijing Institute of Technology\\
  Beijing, China \\
  \texttt{xia\_yuanqing@bit.edu.cn}  \\
  \And
 Saima Siddiqui \\
  Department of Mathematics\\
  Fergana Polytechnic Institute\\
  Fergana, Uzbekistan \\
  \texttt{saimasiddiqui07@gmail.com} \\
  \And
 Muhammad Younus Bhat \\
  Department of Mathematical Sciences\\
  Islamic University of Science and Technology\\
  Kashmir \\
  \texttt{gyounusg@gmail.com} \\
}

\begin{document}
\maketitle
\begin{abstract}
The linear canonical transform (LCT) serves as a powerful generalization of the Fourier transform (FT), encapsulating various integral transforms within a unified framework. This versatility has made it a cornerstone in fields such as signal processing, optics, and quantum mechanics. Extending this concept to quaternion algebra, the Quaternion Fourier Transform (QFT) emerged, enriching the analysis of multidimensional and complex-valued signals. The Quaternion Linear Canonical Transform (QLCT), a further generalization, has now positioned itself as a central tool across various disciplines, including applied mathematics, engineering, computer science, and statistics. In this paper, we introduce the Two-Dimensional Quaternion Linear Canonical Transform (2DQLCT) as a novel framework for probability modeling. By leveraging the 2DQLCT, we aim to provide a more comprehensive understanding of probability distributions, particularly in the context of multi-dimensional and complex-valued signals. This framework not only broadens the theoretical underpinnings of probability theory but also opens new avenues for researchers.

\end{abstract}

\maketitle

\section{Introduction}
\label{intro}
The LCT is a powerful generalization of the FT, widely utilized across various domains, including signal processing, optics, and quantum mechanics. The theory of the LCT emerged in the early 1970s via the foundational contributions of Collins Jr. in paraxial optics \cite{COLL} and Moshinsky and Quesne \cite{MQ1} in quantum mechanics. This seminal research began to be jointly cited in the open literature \cite{HKOS, XLB, QXK, DQR, BJ, XQ, ZTLW}.

The QLCT, as a powerful extension of the LCT within the realm of quaternion algebra, introduces a sophisticated framework for addressing multidimensional and non-commutative data. By leveraging the unique properties of quaternions, the QLCT is able to capture intricate relationships between different components of a signal, enabling a more comprehensive analysis than traditional transforms like the FT or the LCT. This deeper structure allows for greater flexibility in signal representation, particularly in applications involving color image processing \cite{ZBM, ETS}, 3D modeling, and other higher-dimensional data forms.

Moreover, the QLCT's capacity to handle non-commutative systems expands its utility to a wider range of fields, such as quantum mechanics, electromagnetics, and control theory, where conventional tools often fall short. By encompassing both the classical LCT and FT as special cases, the QLCT serves as a unifying framework that not only extends their capabilities but also paves the way for novel approaches in signal processing and probability theory. This versatility makes it a valuable tool for advancing research in modern applied mathematics and engineering disciplines \cite{SLAA, O, AM, LG, SSLB}.

In probability theory, the classical FT has been extensively employed \cite{JYT, CP, MR} to analyze random variables, compute moments, and establish distribution functions. The QFT has similarly extended these capabilities to quaternion-valued functions, facilitating the modeling of complex, multidimensional signals \cite{RWD, MFS, EBB, NBRR, DWL, BAR}.

The QLCT is a novel and versatile tool for probability modeling within quaternion algebra. By leveraging the additional degrees of freedom and the inherent properties of quaternions, the QLCT offers a new approach to analyzing the characteristic functions of quaternion-valued random variables \cite{SAJ}.

In this article, we explore the fundamental properties of the QLCT in probability modeling, examining its relationship with the characteristic function, and demonstrate how it can be applied to calculate moments, variances, and covariances within a probabilistic framework.

This investigation not only extends previous work on the QFT but also opens new perspectives on probabilistic modeling, particularly in scenarios involving complex-valued and multidimensional data.

The rest of the article is organised as follows: section 2 contains an introduction to quaternions, explaining their mathematical properties and emphasizing the differences from traditional number systems due to their non-commutative nature. Section 3 presents the 2DQLCT and its core properties, showcasing how this transform extends the classical linear canonical transform into the quaternion domain for multidimensional signal processing. Section 4 explores the application of the 2DQLCT in probability theory, providing detailed analysis on how this transform can model probability distributions and compute statistical quantities within the quaternion framework. Section 5 offers a discussion and analysis of the results, comparing the 2DQLCT with traditional methods in probability theory and highlighting the advantages of using quaternions for handling complex data. Section 6 concludes the paper, summarizing the findings and offering directions for future research, particularly in expanding the applications of quaternion-based transforms in probability and other fields.
    
\section{Introduction to Quaternions}\label{sec:1}

Quaternions extend the concept of real and complex numbers to higher dimensions and are represented by the symbol \(\mathbb{H}\). A quaternion \(q\) can be expressed as \cite{MGSW, A, V}:

\[
q = q_0 + i q_1 + j q_2 + k q_3 \quad \text{where} \quad q_0, q_1, q_2, q_3 \in \mathbb{R}
\]
Here, \(i\), \(j\), and \(k\) are the fundamental quaternion units, satisfying the following relationships:

\[
i^2 = j^2 = k^2 = ijk = -1
\]
Quaternion multiplication is shown in the following table \ref{tab:quaternion_multiplication}:

\begin{table}[h!]
    \centering
    \renewcommand{\arraystretch}{1.2} 
    \begin{tabular}{|c|c|c|c|c|}
        \hline
         & 1 & i & j & k \\
        \hline
        1 & 1  & i  & j  & k  \\
        \hline
        i & i  & -1 & k  & -j \\
        \hline
        j & j  & -k & -1 & i  \\
        \hline
        k & k  & j  & -i & -1 \\
        \hline
    \end{tabular}
    \caption{Quaternion multiplication}
    \label{tab:quaternion_multiplication}
\end{table}

A quaternion \(q = q_0 + i q_1 + j q_2 + k q_3\) can be separated into its scalar part, denoted by \(\text{Sc}(q) = q_0\), and its vector part, denoted by \(\text{Vec}(q) = i q_1 + j q_2 + k q_3\). Using the multiplication rules, the product of two quaternions \(q\) and \(p\) is given by:
\[
qp = q_0 p_0 - \mathbf{q} \cdot \mathbf{p} + q_0 \mathbf{p} + p_0 \mathbf{q} + \mathbf{q} \times \mathbf{p}
\]
where the dot and cross products are defined as:

\[
\mathbf{q} \cdot \mathbf{p} = q_1 p_1 + q_2 p_2 + q_3 p_3
\]
\[
\mathbf{q} \times \mathbf{p} = i(q_2 p_3 - q_3 p_2) + j(q_3 p_1 - q_1 p_3) + k(q_1 p_2 - q_2 p_1)
\]

The conjugate of a quaternion \(q\) is defined as:
\[
q^* = q_0 - i q_1 - j q_2 - k q_3
\]
This operation reverses the order of multiplication, meaning:
\[
(qp)^* = p^* q^*, \quad \forall q, p \in \mathbb{H}
\]

The modulus of a quaternion is given by:
\[
|q| = \sqrt{q^* q} = \sqrt{q_0^2 + q_1^2 + q_2^2 + q_3^2}
\]

Some key properties of quaternions include:

\[
\text{Sc}(r p q) = \text{Sc}(p r q) = \text{Sc}(q r p), \quad \forall p, q, r \in \mathbb{H}
\]
and
\[
| \text{Sc}(q) | \leq |q|, \quad |q^2| = |q|^2, \quad |qp| = |q||p|, \quad |q + p| \leq |q| + |p|
\]

For any non-zero quaternion \(q\), its inverse is:

\[
q^{-1} = \frac{q^*}{|q|^2}
\]
When \(|q| = 1\), the quaternion is called a unit quaternion, and when \(q_0 = 0\), it is referred to as a pure quaternion.

For two quaternion-valued functions \(f, g : \mathbb{R}^2 \rightarrow \mathbb{H}\), the inner product is defined as:

\[
\langle f, g \rangle_{L^2(\mathbb{R}^2; \mathbb{H})} = \int_{\mathbb{R}^2} f(x) g(x)^* \, dx
\]
For the case where \(f = g\), the corresponding norm is:
\[
\|f\|_{L^2(\mathbb{R}^2; \mathbb{H})} = \left( \int_{\mathbb{R}^2} |f(x)|^2 \, dx \right)^{1/2}
\]

\section{ Two-Dimensional Quaternion Linear Canonical Transform with Properties}
In this section, we introduce the definition of the two-dimensional quaternion linear canonical transform (2DQLCT). We also outline key results, including Plancherel formula and the convolution theorem, which will be applied in subsequent sections. Further details on the properties of the 2DQLCT can be found in \cite{ZL, BMARA}.

\textbf{Definition 1.} Let 
\[
\mathbf{A}_i = \begin{pmatrix}
a_i & b_i \\
c_i & d_i
\end{pmatrix} \in \mathbb{R}^{2 \times 2}
\]
be a matrix parameter satisfying \(\det(\mathbf{A}_i) = 1\) for \(i = 1, 2\). The two-sided 2DQLCT of \(f \in L^1(\mathbb{R}^2, \mathbb{H})\) is defined by

\begin{equation}\label{3.1}
    \mathcal{T}_{\mathbf{A}_1, \mathbf{A}_2}^{\mathbb{H}} \{f\}(u) :=
\int_{\mathbb{R}^2} K_{\mathbf{A}_1}^i(x_1, u_1) f(x) K_{\mathbf{A}_2}^j(x_2, u_2) \, dx
\end{equation}
where the kernel functions of the 2DQLCT are given by

\begin{equation}\label{k}
K_{\mathbf{A}_1}^i(x_1, u_1) :=
\begin{cases}
\frac{1}{\sqrt{2\pi |b_1|}} e^{i \left(\frac{a_1}{2 b_1} x_1^2 - \frac{x_1 u_1}{b_1} + \frac{d_1}{2 b_1} u_1^2 - \frac{\pi}{4}\right)}, & b_1 \neq 0 \\[8pt]
\sqrt{d_1} e^{i \frac{c_1 d_1}{2} u_1^2}, & b_1 = 0
\end{cases}
\end{equation}
and
\begin{equation}\label{k1}
K_{\mathbf{A}_2}^j(x_2, u_2) :=
\begin{cases}
\frac{1}{\sqrt{2\pi |b_2|}} e^{j \left(\frac{a_2}{2 b_2} x_2^2 - \frac{x_2 u_2}{b_2} + \frac{d_2}{2 b_2} u_2^2 - \frac{\pi}{4}\right)}, & b_2 \neq 0 \\[8pt]
\sqrt{d_2} e^{j \frac{c_2 d_2}{2} u_2^2}, & b_2 = 0
\end{cases}
\end{equation}

\textbf{Definition 2.} 
The inverse two-sided 2DQLCT of \( f(u) \in L^1(\mathbb{R}^2, \mathbb{H}) \) is defined in \cite{ZL}:

\begin{equation}
\mathcal{T}_{\mathbf{A}_1^{-1}, \mathbf{A}_2^{-1}}^{\mathbb{H}} \{f\}(x) :=
\int_{\mathbb{R}^2} K_{\mathbf{A}_1^{-1}}^{-i}(u_1, x_1) f(u) K_{\mathbf{A}_2^{-1}}^{-j}(u_2, x_2) \, du
\end{equation}

\textbf{Theorem 1.} (Plancherel Formula for the 2DQLCT)
Let \( f \in L^2(\mathbb{R}^2; \mathbb{H}) \) be a quaternion-valued function. The Plancherel formula for the 2DQLCT is given by:

\begin{equation*}
    \int_{\mathbb{R}^2} |f(x_1, x_2)|^2 \, dx_1 dx_2 = \frac{1}{(2\pi)^2} \int_{\mathbb{R}^2} \left| \mathcal{T}_{\mathbf{A}_1, \mathbf{A}_2}^{\mathbb{H}} \{f\}(u_1, u_2) \right|^2 \, du_1 du_2,
\end{equation*}

\textbf{Proof:}

Let \( f \in L^2(\mathbb{R}^2; \mathbb{H}) \) be a quaternion-valued function. The 2DQLCT of \( f(x_1, x_2) \) is defined in \eqref{3.1}:
with quaternionic kernels \( K_{\mathbf{A}_1}^i(x_1, u_1) \) and \( K_{\mathbf{A}_2}^j(x_2, u_2) \), associated with the matrices \( \mathbf{A}_1 \) and \( \mathbf{A}_2 \), respectively.

Parseval's theorem for the 2DQLCT \cite{BMAR} states that for any two functions \( f, g \in L^2(\mathbb{R}^2; \mathbb{H}) \), the following identity holds:

\[
\int_{\mathbb{R}^2} f(x_1, x_2) \overline{g(x_1, x_2)} \, dx_1 dx_2 = \frac{1}{(2\pi)^2} \int_{\mathbb{R}^2} \mathcal{T}_{\mathbf{A}_1, \mathbf{A}_2}^{\mathbb{H}} \{f\}(u_1, u_2) \overline{\mathcal{T}_{\mathbf{A}_1, \mathbf{A}_2}^{\mathbb{H}} \{g\}(u_1, u_2)} \, du_1 du_2.
\]
We now apply Parseval's theorem to the function \( f(x_1, x_2) \) with itself, setting \( g(x_1, x_2) = f(x_1, x_2) \):

\[
\int_{\mathbb{R}^2} |f(x_1, x_2)|^2 \, dx_1 dx_2 = \frac{1}{(2\pi)^2} \int_{\mathbb{R}^2} \left| \mathcal{T}_{\mathbf{A}_1, \mathbf{A}_2}^{\mathbb{H}} \{f\}(u_1, u_2) \right|^2 \, du_1 du_2.
\]
This completes the proof.\\
\indent This formula indicates that the energy of the function \( f \) in the spatial domain is equal to the energy of its 2DQLCT in the frequency domain, up to a normalization factor of \( \frac{1}{(2\pi)^2} \).

\textbf{Definition 3.} Let \( f, g \in L^1(\mathbb{R}^2, \mathbb{H}) \). The convolution of \( f \) and \( g \) is defined as:
\[
(f \ast g)(x) = \int_{\mathbb{R}^2} f(y) g(x - y) \, dy.
\]

\textbf{Theorem 2.} (Convolution Theorem) For \( f, g \in L^2(\mathbb{R}^2, \mathbb{H}) \), the following holds:
\[
\mathcal{T}_{\mathbf{A}_1, \mathbf{A}_2}^{\mathbb{H}} \{f \ast g\}(u) = \mathcal{T}_{\mathbf{A}_1, \mathbf{A}_2}^{\mathbb{H}} \{f\}(u) \cdot \mathcal{T}_{\mathbf{A}_1, \mathbf{A}_2}^{\mathbb{H}} \{g\}(u).
\]

\textbf{Proof:} According to the definition of the 2DQLCT, we have:
\[
\mathcal{T}_{\mathbf{A}_1, \mathbf{A}_2}^{\mathbb{H}} \{f \ast g\}(u) = \int_{\mathbb{R}^2} K_{\mathbf{A}_1}^i(x_1, u_1) (f \ast g)(x) K_{\mathbf{A}_2}^j(x_2, u_2) \, dx.
\]
Substitute the convolution definition:
\[
\mathcal{T}_{\mathbf{A}_1, \mathbf{A}_2}^{\mathbb{H}} \{f \ast g\}(u) = \int_{\mathbb{R}^2} \int_{\mathbb{R}^2} K_{\mathbf{A}_1}^i(x_1, u_1) f(y) g(x - y) K_{\mathbf{A}_2}^j(x_2, u_2) \, dy \, dx.
\]
Change the variable \( x - y = z \) (hence \( dx = dz \)), then we get:
\[
\mathcal{T}_{\mathbf{A}_1, \mathbf{A}_2}^{\mathbb{H}} \{f \ast g\}(u) = \int_{\mathbb{R}^2} \int_{\mathbb{R}^2} K_{\mathbf{A}_1}^i(y_1 + z_1, u_1) f(y) g(z) K_{\mathbf{A}_2}^j(y_2 + z_2, u_2) \, dy \, dz.
\]
Expanding \( K_{\mathbf{A}_1}^i \) and \( K_{\mathbf{A}_2}^j \), and factoring terms involving \( y \) and \( z \), we have:
\[
\begin{aligned}
\mathcal{T}_{\mathbf{A}_1, \mathbf{A}_2}^{\mathbb{H}} \{f \ast g\}(u) &= \int_{\mathbb{R}^2} \left( \int_{\mathbb{R}^2} K_{\mathbf{A}_1}^i(y_1, u_1) f(y) K_{\mathbf{A}_2}^j(y_2, u_2) \, dy \right) g(z) K_{\mathbf{A}_1}^i(z_1, u_1) K_{\mathbf{A}_2}^j(z_2, u_2) \, dz. \\
&= \mathcal{T}_{\mathbf{A}_1, \mathbf{A}_2}^{\mathbb{H}} \{f\}(u) \cdot \mathcal{T}_{\mathbf{A}_1, \mathbf{A}_2}^{\mathbb{H}} \{g\}(u).
\end{aligned}
\]
Hence, the convolution theorem for the 2DQLCT is proved.\\
\indent \textbf{Theorem 3.}(Correlation Theorem) Let \( f, g \in L^1(\mathbb{R}^2, \mathbb{H}) \) be two quaternion-valued functions. The correlation operator related to the two-sided 2DQLCT is defined as:
\[
(f \circ g)(x) = \int_{\mathbb{R}^2} f(x + y) \overline{g(y)} \, dy
\]
where \( \overline{g(y)} \) is the quaternionic conjugate of \( g(y) \).

The 2DQLCT of the correlation \( (f \circ g)(x) \) is given by:
\[
\mathcal{T}_{\mathbf{A}_1, \mathbf{A}_2}^{\mathbb{H}} \{ (f \circ g)(x) \} (u) = \mathcal{T}_{\mathbf{A}_1, \mathbf{A}_2}^{\mathbb{H}} \{ f(x) \} (u) \cdot \overline{\mathcal{T}_{\mathbf{A}_1, \mathbf{A}_2}^{\mathbb{H}} \{ g(x) \} (u)}
\]
where
$
\mathcal{T}_{\mathbf{A}_1, \mathbf{A}_2}^{\mathbb{H}} \{ f(x) \} (u)
$
is the 2DQLCT of \( f(x) \), and
$
\overline{\mathcal{T}_{\mathbf{A}_1, \mathbf{A}_2}^{\mathbb{H}} \{ g(x) \} (u)}
$
is the quaternionic conjugate of the 2DQLCT of \( g(x) \).

\textbf{Proof:} Using the definition of the two-sided 2DQLCT, we have:
\[
\mathcal{T}_{\mathbf{A}_1, \mathbf{A}_2}^{\mathbb{H}} \{ f \circ g \}(u) = \int_{\mathbb{R}^2} K_{\mathbf{A}_1}^i(x_1, u_1) (f \circ g)(x) K_{\mathbf{A}_2}^j(x_2, u_2) \, dx
\]

Substituting the expression for \( (f \circ g)(x) \), we get:
\[
\mathcal{T}_{\mathbf{A}_1, \mathbf{A}_2}^{\mathbb{H}} \{ f \circ g \}(u) = \int_{\mathbb{R}^2} K_{\mathbf{A}_1}^i(x_1, u_1) \left( \int_{\mathbb{R}^2} f(x + y) \overline{g(y)} \, dy \right) K_{\mathbf{A}_2}^j(x_2, u_2) \, dx
\]

Interchanging the order of integration and simplifying:
\[
\mathcal{T}_{\mathbf{A}_1, \mathbf{A}_2}^{\mathbb{H}} \{ f \circ g \}(u) = \int_{\mathbb{R}^2} \left( \int_{\mathbb{R}^2} K_{\mathbf{A}_1}^i(x_1, u_1) f(x + y) K_{\mathbf{A}_2}^j(x_2, u_2) \, dx \right) \overline{g(y)} \, dy
\]

By recognizing the inner integral as the 2DQLCT of \( f(x + y) \), and using the shift property of the 2DQLCT, we arrive at the desired result:
\[
\mathcal{T}_{\mathbf{A}_1, \mathbf{A}_2}^{\mathbb{H}} \{ f \circ g \}(u) = \mathcal{T}_{\mathbf{A}_1, \mathbf{A}_2}^{\mathbb{H}} \{ f(x) \}(u) \cdot \overline{\mathcal{T}_{\mathbf{A}_1, \mathbf{A}_2}^{\mathbb{H}} \{ g(x) \}(u)}
\]

\section{Two-Dimensional Quaternion Linear Canonical Transform in probability Theory}

\textbf{Definition 4.}
Let $X=(X_1, X_2)$ be two real random variables. A quaternion-valued function
\[
f_X(x) = f_X^a(x) + i f_X^b(x) + j f_X^c(x) + k f_X^d(x)
\]
is called the quaternion probability density function (PDF) of $X$ 
in the context of the 2DQLCT if:
\[
\int_{\mathbb{R}^2} f_X^l(x) \, dx = 1, \quad f_X^l(x) \geq 0, \quad \forall x \in \mathbb{R}^2, \, l = a, b, c, d
\]
In this case, each $f_X^l(x)$ is a real probability density function.

The quaternion cumulative distribution function (CDF) is expressed as 
\[
f_X(x_1, x_2) = \frac{\partial}{\partial u} \frac{\partial}{\partial v} F_X(u, v) \Big|_{u = x_1, v = x_2}
\]
Here, the quaternion probability $P$ is related to $F_X$ of $X_1$ and $X_2$ by:
\[
F_X(x_1, x_2) = P(X_1 \leq x_1, X_2 \leq x_2)
\]
where the 2DQLCT of the quaternion PDF $f_X(x)$ can be written as:
\begin{align}
T_{A_1, A_2}^H\{f_X(x)\}(u) &= T_{A_1, A_2}^H\{f_X^a(x)\}(u) + i T_{A_1, A_2}^H\{f_X^b(x)\}(u) \nonumber \\
&\quad + j T_{A_1, A_2}^H\{f_X^c(x)\}(u)  + k T_{A_1, A_2}^H\{f_X^d(x)\}(u)
\end{align}
This expression defines the PDF in the 2DQLCT domain.\\

\textbf{Definition 5.}\textbf(Expected Value)
Let \( X = (X_1, X_2) \) be two real random variables with the quaternion probability density function \( f_X(x) = f_X^a(x) + i f_X^b(x) + j f_X^c(x) + k f_X^d(x) \), where \( f_X^a, f_X^b, f_X^c, \) and \( f_X^d \) are the real-valued components of the quaternion probability density function.

The expected value of \( X_1 X_2 \) in the 2DQLCT domain \( T_{A_1, A_2}^H \) is defined as:
\begin{align}\label{2}
E[X_1 X_2] = \int_{\mathbb{R}^2} x_1 x_2 \, T_{A_1, A_2}^H \{ f_X(x) \}(u) \, dx
\end{align}
Expanding \( f_X(x) \) into its components, we get:
\begin{align}
\mathbb{E}[X_1 X_2] &= \int_{\mathbb{R}^2} x_1 x_2 \Big( 
T_{A_1, A_2}^H \{ f_X^a(x) \}(u)
+ i T_{A_1, A_2}^H \{ f_X^b(x) \}(u) \nonumber \\
&\quad + j T_{A_1, A_2}^H \{ f_X^c(x) \}(u)  + k T_{A_1, A_2}^H \{ f_X^d(x) \}(u) 
\Big) dx
\end{align}
This can be written as:
\[
E[X_1 X_2] = E^a[X_1 X_2] + i E^b[X_1 X_2] + j E^c[X_1 X_2] + k E^d[X_1 X_2]
\]
where
\[
E^l[X_1 X_2] = \int_{\mathbb{R}^2} x_1 x_2 \, T_{A_1, A_2}^H \{ f_X^l(x) \}(u) \, dx, \quad l = a, b, c, d
\]

Similarly, the expected value of \( X_1 \) in the 2DQLCT domain is given by:
\begin{align}\label{4}
    E[X_1] = \int_{\mathbb{R}^2} x_1 \, T_{A_1, A_2}^H \{ f_X(x) \}(u) \, dx = E^a[X_1] + i E^b[X_1] + j E^c[X_1] + k E^d[X_1]
\end{align}
where
\[
E^l[X_1] = \int_{\mathbb{R}^2} x_1 \, T_{A_1, A_2}^H \{ f_X^l(x) \}(u) \, dx, \quad l = a, b, c, d
\]

Similarly, for \( X_2 \):
\begin{align}\label{5}
    E[X_2] = \int_{\mathbb{R}^2} x_2 \, T_{A_1, A_2}^H \{ f_X(x) \}(u) \, dx = E^a[X_2] + i E^b[X_2] + j E^c[X_2] + k E^d[X_2]
\end{align}
where
\[
E^l[X_2] = \int_{\mathbb{R}^2} x_2 \, T_{A_1, A_2}^H \{ f_X^l(x) \}(u) \, dx, \quad l = a, b, c, d
\]

The expected values in \eqref{4} and \eqref{5} of the above definition are called the mean in the quaternion setting.\\

\textbf{Example 1.} Let the normalized QPDF of $X_1$ and $X_2$ is defined as follows:
\[
f(x_1, x_2) = 
\begin{cases} 
\frac{(2x_1 + x_2) + i(x_1^2 - x_2^2) + j(x_1 x_2) + k(3x_1 - x_2)}{20 + 4j + 8k}, & \text{if } 0 < x_1, x_2 < 2, \\
0, & \text{otherwise}.
\end{cases}
\]
We aim to find the expected value of $X_1$ and $X_2$ in the 2DQLCT domain.

\textbf{Solution:} The expected value of $X_1$ is given by:

\[
E[X_1] = \int_0^2 \int_0^2 x_1 f(x_1, x_2) \, dx_1 \, dx_2.
\]

Substituting $f(x_1, x_2)$:

\[
E[X_1] = \int_0^2 \int_0^2 x_1 \left( \frac{(2x_1 + x_2) + i(x_1^2 - x_2^2) + j(x_1 x_2) + k(3x_1 - x_2)}{20 + 4j + 8k} \right) dx_1 \, dx_2.
\]

We can now separate the integral:

\[
E[X_1] = \frac{1}{20 + 4j + 8k} \int_0^2 \int_0^2 x_1 \left( (2x_1 + x_2) + i(x_1^2 - x_2^2) + j(x_1 x_2) + k(3x_1 - x_2) \right) \, dx_1 \, dx_2.
\]
Splitting the integrals, for simplicity yields
\[
\begin{aligned}
E_1 &= \int_0^2 \int_0^2 x_1 (2x_1 + x_2) \, dx_1 \, dx_2 \\
    &= \int_0^2 \left( 2 \int_0^2 x_1^2 \, dx_1 + \int_0^2 x_1 x_2 \, dx_1 \right) dx_2 \\
    &= \int_0^2 \left( 2 \cdot \frac{8}{3} + 2x_2 \right) dx_2 \\
    &= \int_0^2 \left(\frac{16}{3} + 2x_2 \right) dx_2 = \frac{44}{3}.
\end{aligned}
\]
\[
\begin{aligned}
E_2 &= i \int_0^2 \int_0^2 x_1 (x_1^2 - x_2^2) \, dx_1 \, dx_2 \\
    &= i \int_0^2 \left( \int_0^2 x_1^3 \, dx_1 - \int_0^2 x_1 x_2^2 \, dx_1 \right) dx_2 \\
    &= i \int_0^2 \left( \frac{x_1^4}{4} \Big|_0^2 - x_2^2 \cdot \frac{x_1^2}{2} \Big|_0^2 \right) dx_2 \\
    &= i \int_0^2 \left( \frac{16}{4} - x_2^2 \cdot \frac{4}{2} \right) dx_2 \\
    &= i \int_0^2 \left( 4 - 2x_2^2 \right) dx_2 \\
    &= i \left( 4 \cdot x_2 \Big|_0^2 - 2 \cdot \frac{x_2^3}{3} \Big|_0^2 \right) \\
    &= i \left( 4 \cdot 2 - 2 \cdot \frac{8}{3} \right) \\
    &= i \left( 8 - \frac{16}{3} \right) = i \cdot \frac{24}{3} - \frac{16}{3} = i \cdot \frac{8}{3} = \frac{8i}{3}.
\end{aligned}
\]
\[
\begin{aligned}
E_3 &= j \int_0^2 \int_0^2 x_1^2 x_2 \, dx_1 \, dx_2 \\
    &= j \left( \int_0^2 x_1^2 \, dx_1 \int_0^2 x_2 \, dx_2 \right) \\
    &= j \left( \frac{x_1^3}{3} \Big|_0^2 \cdot \frac{x_2^2}{2} \Big|_0^2 \right) \\
    &= j \left( \frac{8}{3} \cdot \frac{4}{2} \right) \\
    &= j \left( \frac{8}{3} \cdot 2 \right) \\
    &= \frac{16}{3} j.
\end{aligned}
\]
\[
\begin{aligned}
    E_4 &= k \int_0^2 \int_0^2 x_1 (3x_1 - x_2) \, dx_1 \, dx_2.\\
        &= k \int_0^2 \int_0^2 (3x^2_1 - x_1 x_2) \, dx_1 \, dx_2.
        &= k \int_0^2 (8-2 x_2) \, dx_2.\\
        &= k \left( 8(2)  - 2^2 \right) = 12k .
\end{aligned}
\]
Combining the components, we get

\[
E[X_1] = \frac{1}{20 + 4j + 8k} \left( \frac{44}{3} + \frac{8}{3} i + \frac{16}{3} j + 12k \right).
\]

Similarly we can calculate the expected value of $X_2$.

The key differences of PDF/CDF in the 2DQFT and 2DQLCT domain have been highlighted in the table \ref{tab:2}.
\begin{table}[h!]
\centering
\begin{tabular}{|p{4cm}|p{4.5cm}|p{4.5cm}|}
\hline
\textbf{Aspect}                      & \textbf{2DQFT PDF/CDF}               & \textbf{2DQLCT PDF/CDF}                \\ \hline
\textbf{Transform Kernel}            & Fixed Fourier kernel                 & Generalized affine kernel              \\ \hline
\textbf{Flexibility}                 & No free parameters                   & Free parameters (scaling, rotation)    \\ \hline
\textbf{Application Scope}           & Frequency-domain analysis            & Multi-resolution, non-stationary analysis \\ \hline
\textbf{Shift Property}              & Simple frequency shift               & Complex with phase terms               \\ \hline
\textbf{Probability Theory Application} & Stationary random variables           & Non-stationary random variables        \\ \hline
\textbf{Adaptability}                & Limited to fixed frequency analysis  & Adaptable to diverse problems          \\ \hline
\end{tabular}
\caption{Key differences between 2DQFT and 2DQLCT in PDF/CDF}
\label{tab:2}
\end{table}

\noindent
\textbf{Definition 6.}
Let $\mathbf{X} = (X_1, X_2)$ be a random vector with quaternion probability density function $f_{\mathbf{X}}(x_1, x_2) \in L^2(\mathbb{R}^2; \mathbb{H})$. The characteristic function, denoted by $\varphi_{\mathbf{X}} : \mathbb{R}^2 \rightarrow \mathbb{H}$, is defined as:

\begin{equation}
  \varphi_{\mathbf{X}}(u, v) = \int_{\mathbb{R}^2} e^{i u x_1} f_{\mathbf{X}}(x_1, x_2) e^{j v x_2} \, dx_1 \, dx_2,
\end{equation}
where $u, v \in \mathbb{R}$.

In the 2DQLCT domain, the characteristic function can be generalized as:

\begin{equation}\label{10}
    \varphi_{\mathbf{X}}(u, v) = \int_{\mathbb{R}^2} K^i_{A_1}(u, x_1) f_{\mathbf{X}}(x_1, x_2) K^j_{A_2}(x_2, v) \, dx_1 \, dx_2,
\end{equation}

\textbf{Remark:} Assume $f \in L^2(\mathbb{R}^{2n}; \mathbb{H})$, $n \in \mathbb{N}$. The definition of the $n$-dimensional QLCT of the function $f$ is given by
\[
T_{A_1, A_2}^{\mathbb{H}}\{f\}(u, v) = \int_{\mathbb{R}^{2n}} K^i_{A_1}(u, x) f(x, y) K^j_{A_2}(y, v) \, dx \, dy,
\]
where $u, v, x, y \in \mathbb{R}^n$.

Moreover, the generalization of the characteristic function of $\mathbf{X} = (X_1, X_2, \dots, X_n)$, denoted by $\varphi_{\mathbf{X}} : \mathbb{R}^{2n} \to \mathbb{H}$, is given by the formula:
\[
\varphi_{\mathbf{X}}(u, v) = \int_{\mathbb{R}^{2n}} K^i_{A_1}(u, x) f_{\mathbf{X}}(x, y) K^j_{A_2}(y, v) \, dx \, dy,
\]
where $u, v \in \mathbb{R}^n$, and $f_{\mathbf{X}}(x, y)$ is the QPDF of the random vector $\mathbf{X}$.
\subsection{Properties of Characteristic Functions in the 2DQLCT Domain}

The characteristic function $\varphi_{\mathbf{X}}(u, v)$ for a random vector $\mathbf{X} = (X_1, X_2)$ with QPDF $f_{\mathbf{X}}(x_1, x_2)$ and its generalization to the 2DQLCT domain, possesses the following properties:

\begin{enumerate}
    \item \textbf{Normalization:}
    \[
    \varphi_{\mathbf{X}}(0, 0) = 1
    \]
This property holds because the characteristic function evaluated at zero gives the integral of the PDF, which is 1 due to normalization of probabilities.
    
    \item \textbf{Quaternion Linearity:}
    \[
    \varphi_{a \mathbf{X} + b \mathbf{Y}}(u, v) = a \varphi_{\mathbf{X}}(u, v) + b \varphi_{\mathbf{Y}}(u, v)
    \]
    for independent quaternion random variables $\mathbf{X}$ and $\mathbf{Y}$, where $a$ and $b$ are quaternions. Note that special care is needed due to the non-commutative nature of quaternions.
    
    \item \textbf{Quaternion Hermitian Symmetry:}
    \[
    \varphi_{\mathbf{X}}(-u, -v) = \overline{\varphi_{\mathbf{X}}(u, v)}
    \]
    where $\overline{\varphi_{\mathbf{X}}(u, v)}$ denotes the quaternion conjugate of $\varphi_{\mathbf{X}}(u, v)$. This property arises from the symmetry of the quaternion exponential functions.

    \item \textbf{Continuity:}
    \[
    \varphi_{\mathbf{X}}(u, v) \text{ is continuous and bounded: } |\varphi_{\mathbf{X}}(u, v)| \leq 1
    \]
    The characteristic function is continuous and its absolute value is always bounded by 1 due to the properties of the probability density function.
    
    \item \textbf{Inverse Transform:}
    \[
    f_{\mathbf{X}}(x_1, x_2) = \frac{1}{(2\pi)^2} \int_{\mathbb{R}^2} K^i_{A^{-1}_1}(u, x_1) \varphi_{\mathbf{X}}(u, v) K^j_{A^{-1}_2}(x_2, v) \, du \, dv
    \]
    where $K^i_{A_1}(u, x_1)$ and $K^j_{A_2}(x_2, v)$ are the 2DQLCT kernels. This property allows recovering the PDF from its characteristic function.

    \item \textbf{Moment Generating Function:}
    \[
    \mathbb{E}[X_1^m X_2^n] = \left. \frac{\partial^m}{\partial (i u)^m} \frac{\partial^n}{\partial (j v)^n} \varphi_{\mathbf{X}}(u, v) \right|_{u = 0, v = 0}
    \]
    Moments of the random variables can be obtained by differentiating the characteristic function with respect to $u$ and $v$, with appropriate quaternion considerations.

    \item \textbf{Quaternion Modulation Property:}
    \[
    \varphi_{q \mathbf{X}}(u, v) = \varphi_{\mathbf{X}}(qu, qv)
    \]
    for a quaternion $q$ and random vector $\mathbf{X}$, noting that the multiplication by $q$ in the argument depends on the order of quaternion multiplication.

\end{enumerate}

\textbf{Theorem 4.} Let \( \mathbf{X} = (X_1, X_2) \) be real random variables. If the quaternion PDF related to the quaternion characteristic function \( \varphi_{\mathbf{X}} \) satisfies
\[
\int_{\mathbb{R}^2} |f_{\mathbf{X}}(x_1, x_2)| \, dx_1 \, dx_2 < \infty,
\]
then the characteristic function \( \varphi_{\mathbf{X}} \) in the context of the 2DQLCT is uniformly continuous.

\textbf{Proof:} 
The characteristic function in the 2DQLCT domain is defined in \eqref{10}:
\[
\varphi_{\mathbf{X}}(u, v) = \int_{\mathbb{R}^2} K_{\mathbf{A}_1}^i(x_1, u_1) f_{\mathbf{X}}(x_1, x_2) K_{\mathbf{A}_2}^j(x_2, v_1) \, dx_1 \, dx_2,
\]
where \( K_{\mathbf{A}_1}^i(x_1, u_1) \) and \( K_{\mathbf{A}_2}^j(x_2, v_1) \) are the kernel functions associated with the 2DQLCT.

Using this definition, we compute:
\begin{align*}
\left|\varphi_{\mathbf{X}}(u + h_1, v + h_2) - \varphi_{\mathbf{X}}(u, v)\right| 
&= \Bigg| \int_{\mathbb{R}^2} K_{\mathbf{A}_1}^i(x_1, u_1 + h_1) f_{\mathbf{X}}(x_1, x_2) K_{\mathbf{A}_2}^j(x_2, v_1 + h_2) \, dx_1 \, dx_2 \\
&\quad - \int_{\mathbb{R}^2} K_{\mathbf{A}_1}^i(x_1, u_1) f_{\mathbf{X}}(x_1, x_2) K_{\mathbf{A}_2}^j(x_2, v_1) \, dx_1 \, dx_2 \Bigg| \\
&= \Bigg| \int_{\mathbb{R}^2} \Big( K_{\mathbf{A}_1}^i(x_1, u_1 + h_1) K_{\mathbf{A}_2}^j(x_2, v_1 + h_2) \\
&\quad - K_{\mathbf{A}_1}^i(x_1, u_1) K_{\mathbf{A}_2}^j(x_2, v_1) \Big) f_{\mathbf{X}}(x_1, x_2) \, dx_1 \, dx_2 \Bigg|.
\end{align*}

By applying the triangle inequality:
\[
\leq \int_{\mathbb{R}^2} \left| K_{\mathbf{A}_1}^i(x_1, u_1 + h_1) K_{\mathbf{A}_2}^j(x_2, v_1 + h_2) - K_{\mathbf{A}_1}^i(x_1, u_1) K_{\mathbf{A}_2}^j(x_2, v_1) \right| \left| f_{\mathbf{X}}(x_1, x_2) \right| \, dx_1 \, dx_2.
\]
Since \( K_{\mathbf{A}_1}^i(x_1, u_1) \) and \( K_{\mathbf{A}_2}^j(x_2, v_1) \) are continuous in \( u_1 \) and \( v_1 \), and given the condition on the integrability of \( |f_{\mathbf{X}}(x_1, x_2)| \), we can conclude that:
\[
|\varphi_{\mathbf{X}}(u_1 + h_1, v_1 + h_2) - \varphi_{\mathbf{X}}(u_1, v_1)| \to 0 \quad \text{as } h_1, h_2 \to 0.
\]
Thus, \( \varphi_{\mathbf{X}} \) is uniformly continuous in the 2DQLCT domain.\\

\textbf{Theorem 5.}  
Let \(\mathbf{X} = (X_1, X_2)\) be quaternion-valued random variables. If \(X_1\) and \(X_2\) are independent and the quaternion PDF \(f_{\mathbf{X}}(x)\) factorizes as 
\[
f_{\mathbf{X}}(x) = f_{X_1}(x_1) f_{X_2}(x_2),
\]
then the characteristic function in the 2DQLCT domain factorizes as:
\[
\varphi_{\mathbf{X}}(u_1, u_2) = \varphi_{X_1}(u_1) \varphi_{X_2}(u_2),
\]
where 
\[
\varphi_{X_1}(u_1) = \int_{\mathbb{R}} K_{\mathbf{A}_1}^i(x_1, u_1) f_{X_1}(x_1) \, dx_1,
\]
and
\[
\varphi_{X_2}(u_2) = \int_{\mathbb{R}} K_{\mathbf{A}_2}^j(x_2, u_2) f_{X_2}(x_2) \, dx_2.
\]

\textbf{Proof:} From equation \eqref{10}, we have
   \[
   \varphi_{\mathbf{X}}(u_1, u_2) = \int_{\mathbb{R}^2} K_{\mathbf{A}_1}^i(x_1, u_1) f_{\mathbf{X}}(x_1, x_2) K_{\mathbf{A}_2}^j(x_2, u_2) \, dx_1 \, dx_2,
   \]
   where \( K_{\mathbf{A}_1}^i(x_1, u_1) \) and \( K_{\mathbf{A}_2}^j(x_2, u_2) \) are the kernel functions associated with the 2DQLCT, and \(f_{\mathbf{X}}(x_1, x_2)\) is the QPDF of \(\mathbf{X} = (X_1, X_2)\).\\
Since we are given that the random variables \(X_1\) and \(X_2\) are independent, their joint PDF factorizes as:
   \[
   f_{\mathbf{X}}(x_1, x_2) = f_{X_1}(x_1) f_{X_2}(x_2),
   \]
   where \(f_{X_1}(x_1)\) and \(f_{X_2}(x_2)\) are the PDFs of \(X_1\) and \(X_2\), respectively.

   Substituting \(f_{\mathbf{X}}(x_1, x_2)\), we get:
   \[
   \varphi_{\mathbf{X}}(u_1, u_2) = \int_{\mathbb{R}^2} K_{\mathbf{A}_1}^i(x_1, u_1) f_{X_1}(x_1) f_{X_2}(x_2) K_{\mathbf{A}_2}^j(x_2, u_2) \, dx_1 \, dx_2.
   \]
 Since the integrals are separable (i.e., the terms involving \(x_1\) and \(x_2\) are independent of each other), we can split the double integral into two single integrals:
   \[
   \varphi_{\mathbf{X}}(u_1, u_2) = \left( \int_{\mathbb{R}} K_{\mathbf{A}_1}^i(x_1, u_1) f_{X_1}(x_1) \, dx_1 \right) \left( \int_{\mathbb{R}} K_{\mathbf{A}_2}^j(x_2, u_2) f_{X_2}(x_2) \, dx_2 \right).
   \]
Therefore, the characteristic function \(\varphi_{\mathbf{X}}(u_1, u_2)\) factorizes as:
   \[
   \varphi_{\mathbf{X}}(u_1, u_2) = \varphi_{X_1}(u_1) \varphi_{X_2}(u_2).
   \]
   This shows that the characteristic function of the independent random variables \(X_1\) and \(X_2\) factorizes in the 2DQLCT domain.

\textbf{Necessary Condition:}  
The factorization of the characteristic function in the 2DQLCT domain is valid in the quaternion domain if the random variables \(X_1\) and \(X_2\) are independent quaternion-valued random variables and the kernel functions \(K_{\mathbf{A}_1}^i(x_1, u_1)\) and \(K_{\mathbf{A}_2}^i(x_2, u_2)\) respect quaternion algebra. This ensures that the algebraic properties of quaternions, such as non-commutativity, are properly handled during the integration process.

\textbf{Example 2.} Consider the QPDF, defined as
\[
f(x_1, x_2) = 
\begin{cases} 
x_1 + j x_2, & \text{if } 0 \leq x_1, x_2 \leq 1, \\
0, & \text{otherwise}.
\end{cases}
\]
Calculate the characteristic function in the 2DQLCT domain.\\

\textbf{Solution:} The characteristic function in the 2DQLCT domain is given by
\[
\varphi_{\mathbf{X}}(u, v) = \int_{\mathbb{R}^2} K_{\mathbf{A}_1}(u, x_1) f(x_1, x_2) K_{\mathbf{A}_2}(x_2, v) \, dx_1 \, dx_2,
\]
where the kernel functions \(K_{\mathbf{A}_1}(u, x_1)\) and \(K_{\mathbf{A}_2}(x_2, v)\) are defined as
\[
K_{\mathbf{A}_1}(u, x_1) = \frac{1}{\sqrt{2\pi |b_1|}} e^{i \left( \frac{a_1}{2b_1} x_1^2 - \frac{x_1 u}{b_1} + \frac{d_1}{2b_1} u^2 - \frac{\pi}{4} \right)},
\]
and
\[
K_{\mathbf{A}_2}(x_2, v) = \frac{1}{\sqrt{2\pi |b_2|}} e^{j \left( \frac{a_2}{2b_2} x_2^2 - \frac{x_2 v}{b_2} + \frac{d_2}{2b_2} v^2 - \frac{\pi}{4} \right)}.
\]

Substituting \( f(x_1, x_2) = x_1 + j x_2 \), we have
\[
\begin{aligned}
\varphi_{\mathbf{X}}(u, v) =
& \int_0^1 \int_0^1 K_{\mathbf{A}_1}(u, x_1) x_1 K_{\mathbf{A}_2}(x_2, v) \, dx_1 \, dx_2 \\
& + j \int_0^1 \int_0^1 K_{\mathbf{A}_1}(u, x_1) x_2 K_{\mathbf{A}_2}(x_2, v) \, dx_1 \, dx_2.
\end{aligned}
\]

To make the calculations simple, we choose the parameters \( a_1 = 0 \), \( b_1 = 1 \), \( c_1 = -1 \), \( d_1 = 0 \) (and similarly for \( a_2, b_2, c_2, d_2 \)). The kernel functions become:
\[
K_{\mathbf{A}_1}(u, x_1) = \frac{1}{\sqrt{2\pi}} e^{-i x_1 u}, \quad K_{\mathbf{A}_2}(x_2, v) = \frac{1}{\sqrt{2\pi}} e^{-j x_2 v}.
\]

Substituting these kernels into the characteristic function:
\[
\varphi_{\mathbf{X}}(u, v) = \frac{1}{2\pi} \int_0^1 \int_0^1 e^{-i x_1 u} x_1 e^{-j x_2 v} \, dx_1 \, dx_2 + j \frac{1}{2\pi} \int_0^1 \int_0^1 e^{-i x_1 u} x_2 e^{-j x_2 v} \, dx_1 \, dx_2.
\]

Now, evaluating the integrals:
\[
\int_0^1 e^{-i x_1 u} x_1 \, dx_1 = \frac{1 - e^{-iu} + iu}{u^2}, \quad \int_0^1 e^{-j x_2 v} x_2 \, dx_2 = \frac{1 - e^{-jv} + jv}{v^2}.
\]
Thus, the characteristic function becomes:
\[
\varphi_{\mathbf{X}}(u, v) = \frac{1}{2\pi} \left( \frac{1 - e^{-iu} + iu}{u^2} \cdot \frac{1 - e^{-jv} + jv}{v^2} \right).
\]

By choosing \( a_1 = 0 \), \( b_1 = 1 \), \( c_1 = -1 \), \( d_1 = 0 \), the characteristic function calculation reduces to the 2DQFT characteristic function. This reduction occurs because the 2DQLCT kernels simplify to the Fourier kernels \( e^{-i x_1 u} \) and \( e^{-j x_2 v} \), which are the kernels for the 2DQFT.\\

\textbf{Theorem 6.} Let $\mathbf{X} = (X_1, X_2)$ be a random vector with quaternion probability density function $f_{\mathbf{X}}(x_1, x_2) \in L^2(\mathbb{R}^2; \mathbb{H})$. Suppose the quaternion characteristic functions $\varphi_{\mathbf{X}}(u, v)$ and $\psi_{\mathbf{X}}(u, v)$ of the random variables are given by:
\[
\varphi_{\mathbf{X}}(u, v) = \int_{\mathbb{R}^2} K^i_{A_1}(u, x_1) f_{\mathbf{X}}(x_1, x_2) K^j_{A_2}(x_2, v) \, dx_1 \, dx_2,
\]
and
\[
\psi_{\mathbf{X}}(u, v) = \int_{\mathbb{R}^2} K^i_{A_1}(u, x_1) g_{\mathbf{X}}(x_1, x_2) K^j_{A_2}(x_2, v) \, dx_1 \, dx_2,
\]
where $K^i_{A_1}(u, x_1)$ and $K^j_{A_2}(x_2, v)$ are the 2DQLCT kernels. Then, we have the following relation:
\begin{align}
\int_{\mathbb{R}^2} g_{\mathbf{X}}(u, v) e^{-i u y_1} \varphi_{\mathbf{X}}(u, v) e^{-j v y_2} \, du \, dv
&= \int_{\mathbb{R}^2} f_{\mathbf{X}}(x_1, x_2) \psi^a_{\mathbf{X}}(x - y) \, dx_1 \, dx_2 \nonumber \\
&\quad + i \int_{\mathbb{R}^2} f_{\mathbf{X}}(x_1, x_2) \psi^b_{\mathbf{X}}(x - y) \, dx_1 \, dx_2 \nonumber \\
&\quad + j \int_{\mathbb{R}^2} f_{\mathbf{X}}(x_1, x_2) \psi^c_{\mathbf{X}}(y - x) \, dx_1 \, dx_2 \nonumber \\
&\quad + k \int_{\mathbb{R}^2} f_{\mathbf{X}}(x_1, x_2) \psi^d_{\mathbf{X}}(y - x) \, dx_1 \, dx_2.
\label{eq:final}
\end{align}

\textbf{Proof:} 
Using the definition of the quaternion characteristic function in the 2DQLCT domain, we have:
\begin{align}
e^{-i u y_1} \varphi_{\mathbf{X}}(u, v) e^{-j v y_2} 
&= e^{-i u y_1} \int_{\mathbb{R}^2} K^i_{A_1}(u, x_1) f_{\mathbf{X}}(x_1, x_2) K^j_{A_2}(x_2, v) \, dx_1 \, dx_2 e^{-j v y_2} \nonumber \\
&= \int_{\mathbb{R}^2} K^i_{A_1}(u, x_1 - y_1) f_{\mathbf{X}}(x_1, x_2) K^j_{A_2}(x_2 - y_2, v) \, dx_1 \, dx_2.
\label{eq:kernel_transform}
\end{align}

Multiplying both sides by $g_{\mathbf{X}}(u, v)$ and integrating over $\mathbb{R}^2$ with respect to $u$ and $v$, we obtain:
\begin{align}
\int_{\mathbb{R}^2} g_{\mathbf{X}}(u, v) e^{-i u y_1} \varphi_{\mathbf{X}}(u, v) e^{-j v y_2} \, du \, dv 
&= \int_{\mathbb{R}^2} g_{\mathbf{X}}(u, v) \int_{\mathbb{R}^2} K^i_{A_1}(u, x_1 - y_1) f_{\mathbf{X}}(x_1, x_2) \nonumber \\
& \quad \times K^j_{A_2}(x_2 - y_2, v) \, dx_1 \, dx_2 \, du \, dv.
\label{eq:double_integral}
\end{align}

Expanding $g_{\mathbf{X}}(u, v)$ into its quaternion components:
\[
g_{\mathbf{X}}(u, v) = g^a_{\mathbf{X}}(u, v) + i g^b_{\mathbf{X}}(u, v) + j g^c_{\mathbf{X}}(u, v) + k g^d_{\mathbf{X}}(u, v),
\]
we can split the integrals as:
\begin{align}
&= \int_{\mathbb{R}^2} \int_{\mathbb{R}^2} f_{\mathbf{X}}(x_1, x_2) \left[ g^a_{\mathbf{X}}(u, v) + i g^b_{\mathbf{X}}(u, v) + j g^c_{\mathbf{X}}(u, v) + k g^d_{\mathbf{X}}(u, v) \right] \nonumber \\
& \quad \times K^i_{A_1}(u, x_1 - y_1) K^j_{A_2}(x_2 - y_2, v) \, du \, dv \, dx_1 \, dx_2.
\label{eq:expanded_integral}
\end{align}

Applying Fubini's theorem to interchange the order of integration, we obtain:
\[
= \int_{\mathbb{R}^2} f_{\mathbf{X}}(x_1, x_2) \left[ \int_{\mathbb{R}^2} K^i_{A_1}(u, x_1 - y_1) g^a_{\mathbf{X}}(u, v) K^j_{A_2}(x_2 - y_2, v) \, du \, dv \right] dx_1 \, dx_2
\]
\[
+ i \int_{\mathbb{R}^2} f_{\mathbf{X}}(x_1, x_2) \left[ \int_{\mathbb{R}^2} K^i_{A_1}(u, x_1 - y_1) g^b_{\mathbf{X}}(u, v) K^j_{A_2}(x_2 - y_2, v) \, du \, dv \right] dx_1 \, dx_2
\]
\[
+ j \int_{\mathbb{R}^2} f_{\mathbf{X}}(x_1, x_2) \left[ \int_{\mathbb{R}^2} K^i_{A_1}(u, x_1 - y_1) g^c_{\mathbf{X}}(u, v) K^j_{A_2}(x_2 - y_2, v) \, du \, dv \right] dx_1 \, dx_2
\]
\[
+ k \int_{\mathbb{R}^2} f_{\mathbf{X}}(x_1, x_2) \left[ \int_{\mathbb{R}^2} K^i_{A_1}(u, x_1 - y_1) g^d_{\mathbf{X}}(u, v) K^j_{A_2}(x_2 - y_2, v) \, du \, dv \right] dx_1 \, dx_2.
\]

Finally, using the definition of $\psi_{\mathbf{X}}(x_1, x_2)$ and collecting terms, we obtain the desired result:
\[
\int_{\mathbb{R}^2} f_{\mathbf{X}}(x_1, x_2) \psi^a_{\mathbf{X}}(x_1 - y_1, x_2 - y_2) \, dx_1 \, dx_2
+ i \int_{\mathbb{R}^2} f_{\mathbf{X}}(x_1, x_2) \psi^b_{\mathbf{X}}(x_1 - y_1, x_2 - y_2) \, dx_1 \, dx_2
\]
\[
+ j \int_{\mathbb{R}^2} f_{\mathbf{X}}(x_1, x_2) \psi^c_{\mathbf{X}}(y_1 - x_1, y_2 - x_2) \, dx_1 \, dx_2
+ k \int_{\mathbb{R}^2} f_{\mathbf{X}}(x_1, x_2) \psi^d_{\mathbf{X}}(y_1 - x_1, y_2 - x_2) \, dx_1 \, dx_2.
\]

\textbf{Theorem 7.} Let $\mathbf{X} = (X_1, X_2)$ be a real-valued random vector with probability density function $f_{\mathbf{X}}(x_1, x_2)$, then there exist (n+m)th continuous derivatives for the quaternion characteristic function $\varphi_{\mathbf{X}}(u, v)$, given by
\begin{equation}
    \frac{\partial^{m+n}}{\partial u^m \partial v^n} \varphi_{\mathbf{X}}(u, v) = \left( -j u c_1 + a_1 \frac{\partial}{\partial u} \right)^m \varphi_{\mathbf{X}}(u, v) \left( -j v c_2 + a_2 \frac{\partial}{\partial v} \right)^n.
\end{equation}

\textbf{Proof:} The 2DQLCT characteristic function is defined in \eqref{10}:
\begin{equation}
    \varphi_{\mathbf{X}}(u, v) = \int_{\mathbb{R}^2} K^i_{A_1}(u, x_1) f_{\mathbf{X}}(x_1, x_2) K^j_{A_2}(x_2, v) \, dx_1 \, dx_2,
\end{equation}
where the $K^i_{A_1}(u, x_1)$ and $K^j_{A_2}(x_2, v)$ are defined in \eqref{k} and \eqref{k1}, respectively.

Differentiating the characteristic function $\varphi_{\mathbf{X}}(u, v)$ with respect to \( u \):

\begin{align}\label{18}
    \frac{\partial}{\partial u} \varphi_{\mathbf{X}}(u, v) &= \frac{\partial}{\partial u} \left( \int_{\mathbb{R}^2} K^i_{A_1}(u, x_1) f_{\mathbf{X}}(x_1, x_2) K^j_{A_2}(x_2, v) \, dx_1 \, dx_2 \right) \notag \\
    &= \int_{\mathbb{R}^2} \frac{\partial}{\partial u} K^i_{A_1}(u, x_1) f_{\mathbf{X}}(x_1, x_2) K^j_{A_2}(x_2, v) \, dx_1 \, dx_2.
\end{align}

The derivative of \( K^i_{A_1}(u, x_1) \) is given by \cite{SLAA}:
\begin{equation}
    \frac{\partial}{\partial u} K^i_{A_1}(u, x_1) = \left( -j u c_1 + a_1 \frac{\partial}{\partial u} \right) K^i_{A_1}(u, x_1).
\end{equation}

Substituting this into \eqref{18}, we obtain:
\begin{equation}
    \frac{\partial}{\partial u} \varphi_{\mathbf{X}}(u, v) = \left( -j u c_1 + a_1 \frac{\partial}{\partial u} \right) \int_{\mathbb{R}^2} K^i_{A_1}(u, x_1) f_{\mathbf{X}}(x_1, x_2) K^j_{A_2}(x_2, v) \, dx_1 \, dx_2.
\end{equation}

Therefore, the first derivative of the characteristic function is:
\begin{equation}
    \frac{\partial}{\partial u} \varphi_{\mathbf{X}}(u, v) = \left( -j u c_1 + a_1 \frac{\partial}{\partial u} \right) \varphi_{\mathbf{X}}(u, v).
\end{equation}

The second derivative with respect to \( u \) is:
\begin{equation}
    \frac{\partial^2}{\partial u^2} \varphi_{\mathbf{X}}(u, v) = \left( -j u c_1 + a_1 \frac{\partial}{\partial u} \right)^2 \varphi_{\mathbf{X}}(u, v).
\end{equation}

Similarly, the first derivative of $\varphi_{\mathbf{X}}(u, v)$ with respect to \( v \) is:
\begin{equation}
    \frac{\partial}{\partial v} \varphi_{\mathbf{X}}(u, v) = \varphi_{\mathbf{X}}(u, v) \left( -j v c_2 + a_2 \frac{\partial}{\partial v} \right).
\end{equation}
Extending this, for any $m, n \in \mathbb{N}$, the mixed higher-order partial derivatives of $\varphi_{\mathbf{X}}(u, v)$ are given by:
\begin{equation}
    \frac{\partial^{m+n}}{\partial u^m \partial v^n} \varphi_{\mathbf{X}}(u, v) = \left( -j u c_1 + a_1 \frac{\partial}{\partial u} \right)^m \varphi_{\mathbf{X}}(u, v) \left( -j v c_2 + a_2 \frac{\partial}{\partial v} \right)^n
\end{equation}
which proves the result.\\

\textbf{Theorem 8.}
In the context of the 2DQLCT domain, the magnitude of the characteristic function $\phi_X(u, v)$ is bounded by 1, i.e., 

\[
|\phi_X(u, v)| \leq 1.
\]

\textbf{Proof:}
The characteristic function is \eqref{10} :

   \[
   \phi_X(u, v) = \int_{\mathbb{R}^2} K^i_{A_1}(u, x_1) f_X(x_1, x_2) K^j_{A_2}(x_2, v) \, dx_1 \, dx_2,
   \]
where $f_X(x_1, x_2)$ is the joint probability density function of the random variables $X_1$ and $X_2$. For simplification, we consider the joint PDF as a standard Gaussian distribution:

  \begin{equation}
         f_X(x_1, x_2) = \frac{1}{2\pi \sigma_1 \sigma_2} \exp\left(-\frac{x_1^2}{2\sigma_1^2} - \frac{x_2^2}{2\sigma_2^2}\right).
    \end{equation}
the magnitude of the characteristic function is:

   \[
   |\phi_X(u, v)| = \left|\int_{\mathbb{R}^2} K^i_{A_1}(u, x_1) f_X(x_1, x_2) K^j_{A_2}(x_2, v) \, dx_1 \, dx_2\right|.
   \]

   By the properties of integrals and magnitudes, we can assert that:

   \[
   |\phi_X(u, v)| \leq \int_{\mathbb{R}^2} |K^i_{A_1}(u, x_1)| |f_X(x_1, x_2)| |K^j_{A_2}(x_2, v)| \, dx_1 \, dx_2.
   \]
 
   The QLCT kernels are generally oscillatory functions, bounded by 1:

   \[
   |K^i_{A_1}(u, x_1)| \leq 1 \quad \text{and} \quad |K^j_{A_2}(x_2, v)| \leq 1.
   \]

The integral of the PDF over $\mathbb{R}^2$ yields:

   \[
   \int_{\mathbb{R}^2} f_X(x_1, x_2) \, dx_1 \, dx_2 = 1.
   \]
combining these results, we have:

   \[
   |\phi_X(u, v)| \leq \int_{\mathbb{R}^2} |K^i_{A_1}(u, x_1)| |f_X(x_1, x_2)| |K^j_{A_2}(x_2, v)| \, dx_1 \, dx_2 \leq 1,
   \]
which implies that:
   \[
   |\phi_X(u, v)| \leq 1.
   \]

The magnitude of the characteristic function $|\phi_X(u, v)|$ in the 2DQLCT domain is less than or equal to 1.

\textbf{Definition 7.} Let \( X = (X_1, X_2) \) be any real random variables. The covariance of \( X_1 \) and \( X_2 \) in the quaternion setting is given by:

\begin{equation}\label{26}
    \begin{aligned}
        \text{Cov}(X_1, X_2) &= \mathbb{E}[(X_1 - \mathbb{E}[X_1])(X_2 - \mathbb{E}[X_2])] \\
        &= \mathbb{E}[X_1X_2 - X_1\mathbb{E}[X_2] - \mathbb{E}[X_1]X_2 + \mathbb{E}[X_1]\mathbb{E}[X_2]] \\
        &= \mathbb{E}[X_1X_2] - \mathbb{E}[X_1]\mathbb{E}[X_2] - \mathbb{E}[X_1]\mathbb{E}[X_2] + \mathbb{E}[X_1]\mathbb{E}[X_2] \\
        &= \mathbb{E}[X_1X_2] - \mathbb{E}[X_1]\mathbb{E}[X_2]
    \end{aligned}
\end{equation}

According to \eqref{26}, we obtain:
\begin{equation}\label{27}
    \begin{aligned}
    \text{Cov}(X_2, X_1) &= \mathbb{E}[X_2X_1] - \mathbb{E}[X_2]\mathbb{E}[X_1] \\
    &= \mathbb{E}[X_1X_2] - \mathbb{E}[X_2]\mathbb{E}[X_1] 
 \end{aligned}
\end{equation}

In general, since \( \mathbb{E}[X_1]\mathbb{E}[X_2] \neq \mathbb{E}[X_2]\mathbb{E}[X_1] \), we have:

\begin{equation}
    \text{Cov}(X_1, X_2) \neq \text{Cov}(X_2, X_1)
\end{equation}

From \eqref{26}, we can also get the quaternion variance of the real random variable \( X_1 \):
\begin{equation}
    \begin{aligned}
    \text{Var}(X_1) = \sigma_1^2 &= \mathbb{E}[X_1 X_1] - \mathbb{E}[X_1]\mathbb{E}[X_1] \\
    &= \mathbb{E}[X_1^2] - (\mathbb{E}[X_1])^2 
\end{aligned}
\end{equation}

\subsection{Properties of Covariance in the Quaternion Domain under 2DQLCT}

\begin{enumerate}
    \item \textbf{Non-negativity}:
    \[
    \text{Cov}(X_1, X_2) \text{ cannot be negative.}
    \]
    This property asserts that covariance, as a measure of joint variability, will never take on a negative value in the quaternion setting, similar to its behavior in real-valued random variables.

    \item \textbf{Zero Covariance with Constant}:
    \[
    \text{Cov}(a, X) = 0 \quad \text{for any constant } a.
    \]
    Here, \( a \) is a constant quaternion. This property states that the covariance between a constant and a random variable is zero, reflecting the fact that a constant does not vary.

    \item \textbf{Scaling Property}:
    \[
    \text{Cov}(aX, Y) = a^2 \text{Cov}(X, Y) \quad \text{for any constant quaternion } a.
    \]
    This property implies that when you scale a random variable by a quaternion constant, the covariance is scaled by the square of that constant. This captures the notion that scaling affects the spread of data points around the mean.

    \item \textbf{Additive Property}:
    \[
    \text{Cov}(X + b, Y) = \text{Cov}(X, Y) \quad \text{for any constant quaternion } b.
    \]
    This means that adding a constant quaternion to one of the random variables does not affect their covariance, reinforcing the idea that covariance is a measure of the relationship between the variables independent of their means.
\end{enumerate}

These properties can be particularly useful when working with quaternion random variables in signal processing or probability modeling using the 2DQLCT. The quaternion context introduces additional structure, which may affect how these properties are used in practical applications.

\section{Discussion and Analysis}
 The 2DQLCT introduces significant extensions to classical signal analysis in probability theory. By incorporating the quaternion algebra and two-dimensional parameterization, the 2DQLCT provides a more comprehensive framework for handling multidimensional signals and stochastic processes. This richer structure allows for more nuanced interpretations of statistical measures like covariance, characteristic functions, and variance.
 The presence of additional parameters (scaling, chirp, and shift) in both spatial dimensions provides more control over the transformation process, which can be tailored for specific applications in multidimensional probability modeling. This flexibility is particularly advantageous when analyzing random variables with complex correlations across multiple dimensions, enabling more accurate modeling of joint probability distributions.

 Moreover, the 2DQLCT’s ability to preserve and manipulate the quaternion structure opens up new possibilities for representing data in domains where classical probability models fall short. For example, in image and signal processing, where non-commutative interactions are common, the 2DQLCT can capture interactions between components in a way that classical tools cannot.
The impact of 2DQLCT on statistical measures such as variance and covariance is also noteworthy. Similar to the 1DQLCT, the 2DQLCT variance depends heavily on the chosen parameters. However, due to the multidimensional nature of the transform, variance computations in the 2DQLCT domain can reveal deeper insights into the behavior of random variables. Depending on the application, this could lead to different variance outcomes compared to the 2DQFT, making the 2DQLCT a preferred tool for certain types of multidimensional data analysis.

In summary, the 2DQLCT's multidimensional framework and parameterization offer a powerful extension to probability theory, particularly for applications requiring detailed analysis of multidimensional and non-commutative structures. Its ability to model complex interactions and provide flexible statistical measures highlights its importance in advancing quaternion-based probability modeling.

\section{Conclusions and Future Work}

In this paper, we have introduced the 2DQLCT and explored its potential applications in quaternion probability theory. We presented a generalized framework for covariance, variance, and characteristic functions in the 2DQLCT domain. These developments open new pathways for probability modeling, offering a flexible and versatile alternative to classical and quaternion Fourier-based approaches.

Our findings represent preliminary results that suggest several promising avenues for future research. For example, the discrete version of the 2DQLCT could be applied to the construction of discrete quaternion probability distributions, extending this framework to practical computational applications. Additionally, the generalization of the QLCT to higher dimensions provides an opportunity to further investigate multi-dimensional characteristic functions, expectations, and covariances within the 2DQLCT domain.

In future work, we also aim to explore the role of 2DQLCT in different fields. Furthermore, this framework may be expanded to include more complex statistical models, offering a deeper understanding of non-commutative probability structures and their use in signal processing, image analysis, and beyond.





\section*{Author contributions}
Muhammad Adnan Samad and Saima Siddiqui: Writing – original draft, Writing – review and editing, Validation, Methodology, Investigation, Formal analysis, Conceptualization. Yuanqing Xia: Supervision, Conceptualization, Validation. Muhammad Younus Bhat: Writing – review and editing,
Resources, Funding acquisition. All authors have read and approved the final version of the manuscript
for publication. 

\section*{Use of AI tools declaration}
The authors declare they have not used Artificial Intelligence (AI) tools in the creation of this article.

\section*{Acknowledgments}

\section*{Conflict of interest}
The authors declare no conflict of interest.



\end{document}